\documentclass[letterpaper, 10 pt, conference]{ieeeconf}  

\IEEEoverridecommandlockouts                              

\usepackage{amsmath} 
\usepackage{amssymb,bbm}  
\usepackage{algorithm}
\usepackage{bm}
\usepackage{algpseudocode}
\usepackage{dsfont}
\usepackage{hyperref}
\usepackage{graphics}
\usepackage{graphicx}
\usepackage{nicefrac}
\usepackage{xcolor}

\renewcommand{\S}{\mathcal{S}}
\newcommand{\A}{\mathcal{A}}
\newcommand{\bbtheta}{{\bm \Theta}}
\newcommand{\btheta}{{\bm \theta}}
\newcommand{\bQ}{{\bf Q}}
\newcommand{\bbQ}{{\bm {\mathbbm Q}}}
\newcommand{\B}{\mathcal{B}}
\newcommand{\F}{\mathcal{F}}
\newcommand{\bB}{\bm{\mathcal{B}}}
\newcommand{\bPi}{{\Pi}}
\newcommand{\bF}{\bm{\mathcal{F}}}
\newcommand{\M}{\mathcal{M}}
\newcommand{\Mset}{\mathbbm{M}}
\newcommand{\cP}{\mathbbm{P}}
\newcommand{\tpi}{{\tilde{\bPi}}}
\newcommand{\indc}[1]{{\mathbbm{1}_{\left\{#1\right\}}}}
\newcommand{\rpi}{{\Psi}}
\newcommand{\x}[2]{{\bf x}_{#1}^{#2}}
\renewcommand{\a}[2]{{\bf a}_{#1}^{#2} }
\renewcommand{\d}[2]{{\bf d}_{#1}^{#2} }

\newcommand{\lp}[1]{\mbox{LP}(#1)}
\newcommand{\cPii}{{\Pi}_{;i}}
\newcommand{\dPii}{{{\Pi}_{;i}}}
\newcommand{\gi}{{g}_{;i}}
\newcommand{\Psii}{{{\Psi}_{;i}}}
\newtheorem{theorem}{Theorem}
\newtheorem{lemma}{Lemma}
\newcommand{\etapi}{\eta^t_{\bPi,\tpi}}
\newcommand{\eop}{{\hfill $\blacksquare$}}
\newcommand{\mdp}{\mbox{Risk-CMDP}}
\newcommand{\hide}[1]{}

\newcommand{\TR}[2]{#2}  

\title{\LARGE \bf
Fixed-point equations solving Risk-sensitive MDP with constraint
}

\author{Vartika Singh and Veeraruna Kavitha\\
IEOR, IIT Bombay, India
\thanks{*The work of first author is partially supported by Prime minister research fellowship, India}
}

\begin{document}

\maketitle
\thispagestyle{empty}
\pagestyle{empty}

\begin{abstract}
There are no computationally feasible algorithms that provide solutions to the finite horizon Risk-sensitive Constrained Markov Decision Process ($\mdp$) problem, even for problems with moderate horizon. With an aim to design the same, we derive a fixed-point equation such that the optimal policy of $\mdp$ is also a solution.  We further provide two  optimization problems equivalent to the $\mdp$. These formulations are instrumental in designing a global algorithm  that converges to the optimal policy. The proposed  algorithm is based on random restarts and a local improvement step, where the local improvement step utilizes the solution of the derived fixed-point equation; random restarts ensure global optimization.  We also provide numerical examples to illustrate the feasibility of our algorithm for inventory control problem with risk-sensitive cost and constraint. The complexity of the algorithm grows only linearly with the time-horizon.

\end{abstract}


\section{INTRODUCTION}

Classical Markov Decision Process (MDP) problems aim to derive an optimal policy that optimizes the expected combined cost accumulated over time (e.g., \cite{puterman}). Both finite and infinite horizon problems are well studied. These problems can be solved with the help of well-known Dynamic Programming (DP) equations. One can also solve MDP problems using Linear Programming (LP) based approach  (e.g., \cite{puterman}). Further,  the LP formulation facilitates the inclusion of constraints, which is not the case with DP equations.

In many scenarios, it becomes important to consider the variations in accumulated cost for different sample paths and not only the expected value (e.g., \cite{inv,uday}). In this case, the risk-sensitive MDPs (Risk-MDPs) are useful, where higher moments of combined cost are also considered (see \cite{Jaquette73,Jaquette76}). The Risk-MDPs have a sensitivity parameter $\gamma$, called a risk-factor which determines the importance of the higher moments  and the variance. It is well known that, as the risk-factor tends to zero, the value of the Risk-MDP approaches that of the classical MDP (e.g., \cite{Jaquette76,uday_uss}  for discounted cost and \cite{avgg} for average cost problems); the same is also true for the optimal policies. However the policies for higher $\gamma$ can be drastically different.

The Risk-MDP problems can be solved using the corresponding dynamic programming approach (see \cite{Howard}). Recently, authors in \cite{atul} proposed an LP based formulation to solve finite horizon Risk-MDP problems. However, the inclusion of constraints in these problems is not straightforward. The DP equations are again not satisfied; the LP approach in  \cite{atul} can handle Risk-sensitive constrained-MDPs (briefly referred to as Risk-CMDP).  However, in \cite{atul}, the state space of the Risk-MDP is augmented to include the constraints, which causes the state space to grow exponentially over time. Such formulation is impossible to implement for problems with a longer time horizon. To the best of our knowledge, there are no other computationally feasible algorithms in the literature to solve $\mdp$ problem. This calls for an implementable algorithm that can solve $\mdp$ problems, and that does not suffer from the curse of dimensionality. 

For infinite horizon Risk-MDP, in contrast to classical MDPs, stationary policies are not optimal in general (e.g., \cite{Jaquette76,usto,uday_uss}). This makes the derivation of optimal policies even more difficult. In literature, such problems are solved by approximating infinite horizon problems with a finite horizon Risk-MDP; for example, \cite{usto} considers a tail cut-off policy, while \cite{uday_uss} considers a tail replacement policy. These approximations are in line with ultimately stationary policies discussed in \cite{Jaquette76}. 

Including constraints in infinite horizon Risk-MDPs is even more challenging. A recent work  \cite{uday} considers
 infinite horizon $\mdp$ and provides approximate solutions via the solutions of an appropriate finite horizon $\mdp$. This approach provides $\epsilon$-optimal policies which improve as the terminal time of finite horizon $\mdp$ increases to infinity.
 Such approximations require solution of finite-horizon problems with sufficiently large terminal times. This again calls for an implementable algorithm,  that does not suffer from the curse of dimensionality as in \cite{atul}.
 
 One of the main goals of this paper is  an algorithm that solves finite-horizon $\mdp$ problems, whose complexity grows (only) linearly with terminal time.  \textit{The proposed algorithm is based on the derivation of a fixed-point equation which must be satisfied by any optimal policy of the constrained problem.}
 Our contributions are threefold: i)  we derive a fixed-point equation such that the optimal policy of $\mdp$ is also a solution of the fixed-point equation; ii)  we derive two equivalent optimization problems that facilitate the derivation of the optimal policy; and iii) we provide a global iterative algorithm that converges to the optimal policy under certain conditions.
 
To illustrate the feasibility of our approach, we present numerical examples that solve some inventory control problems with constraint. We could easily solve problems with time horizons as long as 1000 decision epochs.

Our approach can easily be extended to the case with  multiple constraints including risk-neutral constraints. We believe we can also extend to a case with a combination of risk-sensitive and risk-neutral objective functions and or constraint functions.

\section{Risk Sensitive MDPs}
Risk-sensitive Markov Decision Process (Risk-MDP) is a sequential decision-making problem that aims to optimize the exponential function of a combination of sequential rewards. In contrast to classical MDPs, it also considers higher moments of the combined cost. As in classical MDP, the risk-sensitive framework consists of time horizon $1,\dots,T$, a finite state space $\S_t$ at time $t$, a finite action space $\A_{t,x}$ at time $t$ and for state $x \in \S_t$, a transition function that determines the probability $p_t(x'|x,a)$ of reaching state $x' \in \S_{t+1} $ based on current state $x \in \S_t$ and action $a \in \A_{t,x}$. In this paper, we consider a finite-horizon problem, i.e., $T<\infty$. A reward $ \tilde{r}_t(x,a,x')$  is achieved at time $t$ when action $a$ chosen in state $x$ results in (next) state $x'$. The aim is to find a policy that optimizes a given objective function constructed from the combined cost. 
Let $\cP(\A)$ be the set of probability distributions on set $\A$.
Any policy consists of decision rules, that prescribe the actions to be taken for a given state and time epoch.  It is represented as ${\bPi} = (d_1,\dots,d_{T-1})$, where $d_t:\S_t \to \cP(\A_{t,\cdot})$ (to be more precise, $x \mapsto \cP(A_{t,x}) $) prescribes the action for $t$-th time-slot; in other words, in this work we restrict our attention to the Markovian policies.

The risk-sensitive objective function under policy $\bPi$ and initial distribution $\alpha_r$ for  the Risk-MDP problem is defined as,

\vspace{-4mm}
{\small \begin{eqnarray}
    J_r(\bPi, \alpha_r) &=& E^\bPi_{\alpha_r}\left[e^{ \sum_{t=1}^{T-1}r_t(X_t,A_t,X_{t+1}) +  r_T(X_T)}\right], \mbox{\normalsize with, }  \nonumber \\
    r_t(x,a ,x') &:=& \gamma \beta^t \tilde{r}_t(x, a, x') \label{eqn_J_r} \mbox{ \normalsize and }r_T(x) = \gamma\beta^T  \tilde{r}_T(x).
\end{eqnarray}}%
Here, $\gamma$ is the risk-factor, $\beta \in (0,1)$ is the discount factor,  $(X_t,A_t)_{t\le T-1}$ is the stochastic state-action trajectory that evolves under policy $\bPi$, and $E^\bPi_{\alpha_r}$ represents the expectation under policy $\bPi$ and initial distribution $\alpha_r$. A higher $|\gamma|$ indicates more importance to the higher moments of the combined cost, while with $\gamma \to 0$, one can approach the classical MDP problem. The aim in Risk-MDP problems is to find an optimal policy, i.e., a policy $\bPi^*$ that satisfies, 
\begin{equation}\label{eqn_obj_unconstraint}
 \bPi^* \in \arg\sup_{\bPi\in \Gamma } J_r(\bPi, \alpha_r),
\end{equation}%
where $\Gamma$ is the set of Markovian randomized policies.

The value function $u_t(x)$ at time $t$ and in state $x\in\S_t$ is defined to be (see \cite{Howard,atul}),

\vspace{-3mm}
{\small
\begin{eqnarray}\label{eqn_value_unconstraint}
u_t(x) &=&  \sup_{\bPi} J_r(t,\bPi,x),\mbox{ \normalsize where }\\
J_r(t,\bPi,x)\hspace{-3mm} &:=& \hspace{-3mm} E^{\bPi}_{\alpha_r}\left[e^{ \sum_{\tau=t}^{T-1}r_\tau(X_\tau,A_\tau,X_{\tau+1}) +  r_T(X_T) }|X_t=x\right].\nonumber
\end{eqnarray}}%
By strong Markov property applicable under Markovian policies, the above quantity depends only on sub-policy (of $\bPi$) from $t$ onwards. Observe here that, the optimal policy $\bPi^*$ is the one that achieves the value function $u_1(x)$ for all $x\in \S_1$.
The well-known DP equations to solve the Risk-MDPs are  as follows (for $x\in \S_t$),

\vspace{-3mm}
{\small
\begin{eqnarray}\label{eqn_dp_uncst}
u_T(x) &=& e^{ r_T(x)} \mbox{ \normalsize and for } t\le T-1,\\ 
u_t(x) &=& \max_a \left\{  \sum_{x'}e^{r_t(x,a,x')} p_t(x'|x,a) u_{t+1}(x')\right\}. \nonumber
\end{eqnarray}}%

The fixed point dynamic programming equations  \eqref{eqn_dp_uncst} facilitate the derivation of optimal policy for unconstrained problems. However, such equations are not known for constrained problems, which we introduce in the immediate following. One of the aims of this paper is to  derive an appropriate fixed-point equation that solves the $\mdp$.

\subsection{Risk-sensitive constrained  MDP ($\mdp$)}
We now consider a constraint in  the Risk-MDP problem defined in \eqref{eqn_obj_unconstraint}. Here, at time $t$, an immediate constraint-cost ${\tilde c}_t(x,a,x')$ is incurred along with the reward ${\tilde r}_t(x,a,x')$, when  action $a$  chosen in state $x$ results in next state $x'$. The aim is to keep the expected exponential of the combined  constraint-cost below a certain bound. To keep it general, let $\beta_c$ and $\gamma_c$  be the discount and risk-factors corresponding to the constraint; the factors $\gamma_c$, $\beta_c$ can be different from respective factors $\gamma, \beta$ corresponding to  ${\tilde r}_t$. Thus a policy $\bPi$ is feasible if it satisfies the following constraint,

\vspace{-3mm}
{\small
\begin{eqnarray}\label{eqn_constraint}
J_c(\bPi,\alpha_c) \hspace{-2mm} &=& \hspace{-2mm}\hspace{-2mm} E_{\alpha_c}^\bPi\left[e^{ \sum_{t=1}^{T-1}{c}_t(X_t,A_t,X_{t+1}) +  c_T(X_T)}\right]\le B, \mbox{\normalsize with,}  \nonumber\\
c_t(x,a,x') \hspace{-2mm}&:=& \hspace{-2mm}\gamma_c \beta_c^t \tilde{c}_t(x,a,x'),\  \mbox{\normalsize and }  c_T(x ) \:= \ \gamma_c \beta_c^T  \tilde{c}_T(x),
\end{eqnarray}}%
where $\alpha_c$ is the initial distribution and can be different from $\alpha_r$, the initial distribution of state corresponding to ${\tilde r}_t$.  Thus, the overall problem is,
\begin{eqnarray}\label{eqn_original_problem}
\sup_{\bPi\in \Gamma} && J_r(\bPi,\alpha_r),\\
\mbox{subject to} &&  J_c(\bPi,\alpha_c) \le B. \nonumber
\end{eqnarray}%
Let $\Gamma_c := \{ \bPi \in \Gamma: J_c(\bPi,\alpha_c) \le B \}$ represent the corresponding feasible region. 
\textit{Throughout we assume the existence of a solution for \eqref{eqn_original_problem},} and the aim is to design an algorithm that obtains the same. Towards this, as a first step
we derive a fixed-point equation, in the next section  whose solution optimizes  \eqref{eqn_original_problem}. 


\section{Fixed-point equation}
We begin this section with a few definitions. For ease of notation, we let  $m_t$ represent the immediate reward $r_t$ or constraint-cost $c_t$ function at time $t$ depending upon the choice  $m\in\{r,c\}$. Thus $m_t(x,a,x')$  is the reward/constraint-cost  function at time $t$, when action $a$  chosen in state $x$ results in   $x'$, the new state.

\subsection{Forward factors} 
  For any policy $\bPi$,  time $t$, reward/constraint-cost $m$ and state  $x_t$, define the forward factors $\theta^\bPi_{m,t}(x_t)$ as below:

 \vspace{-2mm}
 {\small
\begin{eqnarray*}
\theta^\bPi_{m,t} (x_t)&:=& E^\bPi_{\alpha_m}\left[e^{\sum_{k=1}^{t-1}m_k(X_k,A_k,X_{k+1})} \indc{X_t =x_{t}}\right],\\
&& \hspace{-22mm}
 = \hspace{-3mm}\sum_{ \x{1}{t-1},\a{1}{t-1}}\hspace{-2mm} \alpha_m(x_1) \hspace{-2mm} \prod_{k\le t-1}  d_k(x_k,a_k) p_k(x_{k+1}| x_k, a_k) e^ { m_k(x_k,a_k,x_{k+1})},
\end{eqnarray*}}%
where $\x{1}{t-1}:= (x_1,\dots,x_{t-1})$ is a vector with  each  $x_k \in {\S}_k$ 
and $\a{1}{t-1}:= (a_1,\dots,a_{t-1})$ is a vector with  each  $a_k \in {\A}_{k,x_k}$.  These factors represent the expected reward/constraint-cost accumulated  under policy $\bPi$ till time $t-1$, and the probability that $X_t=x_{t}$. For any $t$, $m$, it is easy to verify that the forward factors satisfy the following recursive equations for   any $x' \in \S_{t+1}$:
\begin{eqnarray}\label{eqn_mu}
&&\hspace{-10mm}\theta^\bPi_{m,1} (x)= \alpha_m(x), \mbox{ for all} \ x \in \S_1\mbox{ and for } t\ge 1, \\
&&\hspace{-10mm}\theta^\bPi_{m,t+1} (x')= \sum_{x,a}\theta^\bPi_{m,t} (x) d_t(x,a) p_t(x'| x, a) e^{m_t(x,a,x')}.\nonumber
    \end{eqnarray}%
Define the corresponding vectors,  $\bbtheta^\bPi_m := \{ \btheta_{m,t}^\bPi \}_t $  where vector for time $t$,  
$\btheta_{m,t}^\bPi :=\{ (\theta^\bPi_{m,t} (x)) : x \in \S_t\}
$.%
\hide{\color{red}For any $\bPi$, observe that the above defined vectors $\bbtheta_m^\bPi$ satisfy the following functional equation for any $t \ge 1$,
\begin{eqnarray*}
\theta^\bPi_{m,t+1} (x) =\F^\bPi_{m,t+1} (\bbtheta^\bPi_{m};x) \mbox{  for all }  x \in \S_{t+1},
\end{eqnarray*}%
where  the operators $\bF^\bPi_m:=(\F^\bPi_{m,1},\dots,\F^\bPi_{m,T})$ are defined component-wise as below (for any $x' \in \S_t$):
\begin{eqnarray*}
&& \hspace{-10mm}\F^\bPi_{m,1} (\bbtheta;x'):= \theta_{1,x'} = \alpha_m(x'), \mbox{ and for } t\ge 2, \\
&& \hspace{-10mm}\F^\bPi_{m,t+1} (\bbtheta;x'):= \sum_{x,a} \theta_{t,x}  d_t(x,a) p(x'| x, a) e^{ m_t(x,a,x')},
\end{eqnarray*}}

%
%

\subsection{Backward factors} For any time $t$, reward/constraint-cost $m$ and policy  $\bPi$, define the backward factors, $Q^\bPi_{m,t}(x,a)$ for  $x \in \S_t$ and $a \in \A_{t,x}$ as follows,

\vspace{-4mm}
{\small\begin{eqnarray}\label{eqn_q_factor_exp}
    Q^\bPi_{m,t}(x,a) \hspace{-15mm}\\
    &=& E_{\alpha_m}^\bPi\left[e^{\sum_{\tau=t}^{T-1}m_\tau(X_\tau,A_\tau,X_{\tau+1}) + m_T(X_T)}| X_t=x,A_t=a \right], \nonumber
\end{eqnarray}}%
and observe that these   factors satisfy (e.g., see \cite{Howard,atul}),

\vspace{-3mm}
{\small\begin{eqnarray}\label{eqn_Q_factors}
&& \hspace{-1cm} Q^\bPi_{m,T} (x,a) = e^{ m_T(x)}, \mbox{\normalsize and for } t\le T-1,\\
&& \hspace{-1cm} Q^\bPi_{m,t} (x, a) = \sum_{x',a'}  e^{m_t(x,a,x')} p_t(x'|x,a)  d_t(x',a')Q^\bPi_{m,t+1}(x',a'). \nonumber
\end{eqnarray}}%
These backward factors $\{ Q^\bPi_{m,t}(x,a)\}$ represent  the  well known Q-factors for Risk-MDPs, which equal the  `cost-to-go'  from time $t$ onwards, given   $X_t=x$, $A_t=a$ and when policy $\bPi$ is used from $t+1$ onwards. Observe that \eqref{eqn_Q_factors} is similar to the \textit{policy-evaluation step} for classical MDPs (e.g., \cite{puterman}). %
Also define the vectors,
\begin{eqnarray*}
    \bbQ^\bPi_m &:=& \{ \bQ_{m,t}^\bPi \}_t  \mbox{ where, }\\
    \bQ_{m,t}^\bPi &:=&\{ (Q^\bPi_{m,t} (x,a)) : x \in \S_t, a \in \A_{t,x}\}.
\end{eqnarray*}
 \hide{\color{red}For any $\bPi$, the backward factors satisfy the following functional equations, 
\begin{eqnarray*}
Q^\bPi_{m,t} (x,a) =\B^\bPi_{m,t} (\bbQ^\bPi_m;x,a),
\end{eqnarray*}%
where the operators $\bB^\bPi:=(\B^\bPi_{m,1},\dots,\B^\bPi_{m,T})$ are defined component-wise as below (for any $(x,a)$):

\vspace{-3mm}
{\small \begin{eqnarray*}
&& \hspace{-10mm}\B^\bPi_{m,T} (\bbQ;x,a):= e^{m_T(s)}, \mbox{ and for } t\le T-1, \\
&& \hspace{-10mm}\B^\bPi_{m,t} (\bbQ;x,a):=  e^{m_t(x,a,x')} \sum_{s',a'} p(s'|x,a)  d_t(s',a')Q_{t+1}(s',a').
\end{eqnarray*}}}%
\subsection{Linear Program and Fixed-point equation}
A Linear program (LP) is an important ingredient of our proposed fixed-point equation. We now discuss the corresponding objective function.
Let $\tpi=(\tilde{d}_1,\cdots,\tilde{d}_{T-1})$ and $\bPi = ({d}_1,\cdots,{d}_{T-1})$ be any two  policies, and  define the following function indexed by $t\le T-1$,

\vspace{-3mm}
{\small
\begin{equation}\label{eqn_f_lin}
    f_t(\tpi, \bbtheta^\bPi_m, \bbQ^\bPi_m) : =  \sum_{x,a}  \theta^\bPi_{m,t}(x) \tilde{d}_t(x, a) Q^\bPi_{m,t}(x,a).
\end{equation}}%
It is easy to verify that, $\tpi \mapsto f_t(\tpi,\bbtheta,\bbQ)$ is a linear function once $(\bbtheta,\bbQ)$ are fixed (for any $t$).

Let $\etapi:=(\d{1}{t-1},\tilde{d}_t,
\d{t+1}{T-1})$ be  a policy that differs from $\bPi$ only at $t$-th epoch, at which the decision is taken according to policy $\tpi$. We now show that the linear function in \eqref{eqn_f_lin} equals the risk-sensitive cost $J_m$ under policy $\etapi$ and initial distribution $\alpha_m$ (see \eqref{eqn_J_r} or \eqref{eqn_constraint}).
%
\begin{lemma}\label{lem_exp_rew_pi_tilde}{\it For any given pair  of policies $\bPi$ and $\tpi$, and time $t\le T-1$ we have,
\begin{equation*}
  f_t(\tpi,\bbtheta^\bPi_m, \bbQ^\bPi_m ) = J_m(\etapi, \alpha_m). 
\end{equation*}}
\end{lemma}
\noindent {\bf Proof} is in \TR{\cite{TR}}{Appendix}. \eop

Thus one can capture the value of
risk-sensitive objective/constraint function for all the policies that   deviate from $\bPi$ at $t$ using the $f_t$ function.  This observation is crucial in deriving the required fixed-point equation. Towards this, we define one LP  for each policy  $\bPi$ as below: 
\begin{eqnarray*}\label{eqn_LP_formulation}
\lp{\bPi}:&& \max_{\tpi} \sum_{t=1}^{T-1} f_t(\tpi, \bbtheta_r^\bPi, \bbQ^\bPi_r),\\
&& \hspace{-6mm}\mbox{subject to } f_t(\tpi, \bbtheta_c^\bPi, \bbQ^\bPi_c)\le B \mbox{ for all } 1 \le t < T.
\end{eqnarray*}%
When $\bPi=\tpi$, clearly $\etapi=\bPi$, and then from  Lemma \ref{lem_exp_rew_pi_tilde}, 
the  function $f_t$  equals the risk-sensitive objective/constraint function: 
\begin{equation}\label{eqn_obj_val}
    f_t(\bPi,\bbtheta^\bPi_m, \bbQ^\bPi_m ) = J_m(\bPi, \alpha_m) \mbox{ for all }t.
\end{equation}%
Thus if $\bPi$ is feasible for $\mdp$ \eqref{eqn_original_problem}, then it is also feasible for $\lp{\bPi}$ and vice-versa.

Now consider $\bPi^*$, a solution of $\mdp$ \eqref{eqn_original_problem}, and consider $\lp{\bPi^*}$. One can anticipate that $\bPi^*$ solves the $\lp{\bPi^*}$. This is indeed true and provides the required fixed-point equation as shown in the following, which is proved with the help of Lemma \ref{lem_exp_rew_pi_tilde}.
 
 %
%
\begin{theorem}[Necessary condition]\label{thm_nec_cndi}
{\it Let $\M(\bPi)$  be  the solution set of $\lp{\bPi}$, for any $\bPi$. Then, any optimal policy $\bPi^*$ of $\mdp$ \eqref{eqn_original_problem} satisfies the fixed-point equation, }
\begin{equation}\label{eqn_the_fixed_pt}
    \bPi \in \M (\bPi).
\end{equation} 
\end{theorem}
\noindent {\bf Proof} is in \TR{\cite{TR}}{Appendix}. \eop

The above theorem provides a necessary condition to be satisfied by an optimal policy $\bPi^*$.
In the next, we provide two optimization problems  that  are equivalent to the $\mdp$ problem \eqref{eqn_original_problem}, which lead to the required solution/algorithm.

\subsection{Solutions of $\mdp$}
The solution of  fixed-point equation \eqref{eqn_the_fixed_pt} 
is guaranteed to exist, once $\mdp$ \eqref{eqn_original_problem} has a solution; however, it may not be unique. 
Let $\Mset:= \{\bPi : \bPi \in \M(\bPi)\}$ be the set of all such possible fixed points.   We now have our main result,
\begin{theorem}\label{thm_eqvlnt}
{\it 
The solution of the $\mdp$ problem \eqref{eqn_original_problem} is obtained by solving any of the  following two (global fixed point and global optimization) problems,}
\begin{align*}
     \mbox{\bf GF} :       &   \sup_{ \bPi \in \Mset} \sum_{t=1}^{T-1} f_t(\bPi, \bbtheta_r^\bPi, \bbQ_r^\bPi). \\ 
     \mbox{\bf GO}  :  &  
            \sup_{ \bPi \in \Gamma}  \sum_{t=1}^{T-1} f_t(\bPi, \bbtheta_r^\bPi, \bbQ_r^\bPi ), & \\  
            \vspace{2mm}
            && \hspace{-30mm}  \mbox{ \normalfont s. t. } f_t(\bPi, \bbtheta_c^\bPi, \bbQ_c^\bPi) \le B \mbox{ for all }t\le T-1.
\end{align*}
\end{theorem}
\noindent {\bf Proof} is in \TR{\cite{TR}}{Appendix}. \eop

The {\bf GO} problem is just a restatement of the $\mdp$ problem \eqref{eqn_original_problem}, while the equivalence of {\bf GF} problem is proved using the fixed points of Theorem \ref{thm_nec_cndi}.

From Theorem \ref{thm_nec_cndi},   the optimal policy satisfies the fixed-point equation $\bPi \in \M(\bPi)$, and is an element of $\Mset$.  
Hence an iterative algorithm that converges to the fixed points in $\Mset$ can be utilized to derive the solution of $\mdp$ \eqref{eqn_original_problem}. 
However, as seen from \textbf{GF} problem of Theorem \ref{thm_eqvlnt}, one needs to converge towards the best   among the set of fixed points in $\Mset$.  It is important to observe here that both \textbf{GF} and \textbf{GO} formulations maximize the same objective function. Thus a global optimization problem constructed using \textbf{GO} formulation can be of help. 

\section{Algorithm}

 By \textbf{GO} problem of Theorem \ref{thm_eqvlnt}, the $\mdp$ is converted into a  constrained global optimization problem. This optimization problem can be solved using any random search method (e.g., random restarts \cite{global}, simulated annealing \cite{sim_annealing} etc.), provided it satisfies some
 regularity conditions. A global algorithm with random restarts (e.g.,\cite{global}) has two types of update steps: i) a purely random search step -- a random new point is chosen  from   the feasible region at such iterative step, and, ii) a local improvement step -- an appropriate algorithm (e.g. gradient descent) improves the previous update using the new observations. In any iteration $k$, the random restart step is chosen with a certain probability $p_k$, where $p_k$ diminishes with $k$.

We design the local improvement step for the \textbf{GO} problem using the \textbf{GF} problem; combining it with random restarts, we construct a global algorithm that reaches the best fixed point in $\Mset$.
%
%
%
%

\subsection{Local Improvement algorithm}
To begin with, we propose a local improvement step and derive its analysis when it runs continually unperturbed by the random restarts. 
The aim in the local improvement step is to  converge to a fixed point in $\Mset$.  Towards this we propose an iterative algorithm, where the update  for any $k$ is given by the following:
\begin{eqnarray}\label{eqn_local_update}
    \bPi_{k+1}& = &  \bPi_k+ \epsilon_k (\rpi_{k}(\bPi_k) -\bPi_k),  \mbox{ where, }\\
    \rpi_k(\bPi_k) &\in& \M(\bPi_k) \mbox{ and, }\nonumber\\
    \epsilon_k \in (0,1), && \sum_{k=1}^\infty \epsilon_k = \infty, \ \sum_{k=1}^\infty \epsilon_k^2 < \infty. \nonumber
\end{eqnarray}In the above equation, $\rpi_{k}(\bPi_k)$ is chosen randomly from $\M(\bPi_k)$, the solution set of $\lp{\bPi_k}$.  We also allow the LP solver (at $\bPi_k$) to return   a random approximate solution $\rpi_k(\bPi_k)$ whose expected value equals one of the   solutions in $\M(\bPi_k)$. All we require is that the solver solution satisfies the following:
$$
E[ \rpi_k (\bPi_k) | \bPi_k]  \in \M (\bPi_k) \mbox{ almost surely}.
$$

We prove that the process $\{\bPi_k\}_k$ converges to a locally asymptotically stable set $\mathbb{A}$ (set of attractors\footnote{We say an equilibrium point of an ODE is  attractor if it is locally asymptotically stable in the sense of Lyapunov.}) of the following ordinary differential equation (ODE) in Theorem \ref{thm_cnvg_ode} 
(given below) under certain assumptions.
\begin{equation}\label{eqn_ode}
\dot{\bPi} = g(\bPi) =  E[\rpi(\bPi)] - \bPi.
\end{equation}%
We prove this theorem using stochastic approximation-based tools (\cite{kushner2003stochastic}).
Towards this, 
we assume the following:

\noindent $\bullet$ {\bf B} The solution $\bPi^*$ of $\mdp$ \eqref{eqn_original_problem} is unique solution of the corresponding $\lp{\bPi^*}$. The function $g(\cdot)$ is measurable,  and   the optimal policy $\bPi^*$ is an attractor for ODE \eqref{eqn_ode}.

Observe that any equilibrium point $\bPi$ of the above ODE corresponds   to a fixed point in $\Mset$ if $\M(\bPi)=\{\bPi\}$ (note $\M(\bPi)$ is   also a singleton, and hence  $E[\Psi(\bPi)]$ in \eqref{eqn_local_update} equals  $\bPi$).  The above assumption hence requires that $\M(\bPi^*) = \{\bPi^*\}$ at $\bPi^*$, the solution of $\mdp$ \eqref{eqn_original_problem}. Our algorithm may not work without this assumption. 

\begin{theorem}\label{thm_cnvg_ode}
{\it Assume {\bf B}.
Let $\mathbb{A}$ be the set of attractors. Suppose
that $\{\bPi_k\}$ defined in \eqref{eqn_local_update}  visits a compact set in the domain of attraction of $\mathbb{A}$ infinitely
often with probability $\rho$. Then $\bPi_k \to \mathbb{A}$ with probability at least $\rho$.} 
\end{theorem}
\noindent {\bf Proof} is in \TR{\cite{TR}}{Appendix}. \eop

Thus when the local improvement algorithm is not perturbed by random restart steps, the algorithm in \eqref{eqn_local_update} converges to one of the fixed points in $\Mset$; only equilibrium points of \eqref{eqn_ode} can be the attractors in $\mathbb{A}$, and any equilibrium point is a fixed point\footnote{This may not be true if  there is a $\tpi$ such that by randomization, $E[\rpi(\tpi)|\tpi]=\tpi$; we assume this is not the case in the current paper.} in $\Mset$. However this requires that  the algorithm visits a neighborhood of the fixed points (in $\Mset$) infinitely often. We will next analyze the global algorithm and show that the global algorithm visits any such neighborhood infinitely often. Of course, one requires technical proof to connect the two results, which can be a part of future work. For now, we proceed with the analysis of the global algorithm.

\subsection{Global algorithm}
Above theorem shows that the iterates in \eqref{eqn_local_update} converge to the attractors of ODE \eqref{eqn_ode} under certain conditions. These attractors in turn correspond to fixed points in $\Mset$. Recall that set $\Mset$ may contain many other fixed points, which are not optimal. Our aim is to avoid these local points, and arrive at the global optimizer.


With an aim to derive the global optimizer of \textbf{GO} problem, or equivalently the optimal policy for $\mdp$ \eqref{eqn_original_problem}, we propose an algorithm in \ref{alg_random_search},  namely (GRC) Global $\mdp$ algorithm. This algorithm uses the random restarts technique to get closer to the optimal policy along with the local improvement step of \eqref{eqn_local_update} to converge to the best fixed point in $\Mset$. At   iterate $k$, a random policy is chosen  independently from the space of policies $\Gamma$ according to   Uniform distribution  ${\cal U}$,  with probability $p_k$. With the remaining probability, we perform the local improvement update step over the current policy using \eqref{eqn_local_update}. At every iterate, the best policy seen so far is stored.

\begin{algorithm}
\caption{Global $\mdp$ algorithm (GRC)}\label{alg_random_search}
Initialize $\bPi_0$ randomly,  set $J^*_r= -\infty$, $\hat{\bPi}^*=\bPi_0$, choose a constant $w$

\textbf{For} $k=1,2,\dots$
\begin{algorithmic}
\State \vspace{-5mm}\begin{eqnarray*}
\bPi_k \hspace{-1mm} \gets \hspace{-1mm} \left\{ \hspace{-1mm}
\begin{array}{ll}
     \mbox{random policy chosen according to } {\cal U}&  \mbox{w.p. } p_k = \frac{w}{k} \\
     \mbox{Local improvement }  (\bPi_{k-1}) \mbox{ of } \eqref{eqn_local_update}&  \mbox{w.p. } 1-p_k \\ 
\end{array}\right.
\end{eqnarray*}
\State Calculate $J_c(\bPi_k, \alpha_c)$ using \eqref{eqn_obj_val}

\If{$J_c(\bPi_k, \alpha_c) \le B$}
\State Calculate $J_r(\bPi_k, \alpha_r)$ using \eqref{eqn_obj_val}
\vspace{-4mm}
\State \If{$J_r(\bPi_k, \alpha_r) \le J^*_r$}
\State \begin{equation} \label{eqn_pi_hat}
   \hat{\bPi}^*\gets \bPi_k \hspace{4.5cm}
\end{equation}
\State $J^*_r \gets J_r(\bPi_k, \alpha_r)$
\EndIf
\Else 
\State $\Pi_k \gets $ random policy  according to ${\cal U}$
\EndIf
\end{algorithmic}
\end{algorithm}

The Global algorithm  \ref{alg_random_search} converges to the optimal policy as shown in theorem \ref{thm_conv_randm_search}, when the support of restart distribution ${\cal U}$ is entire $\Gamma$, the space of Markovian policies. Also, recall 
 $\Gamma_c := \{ \bPi \in \Gamma: J_c(\bPi,\alpha_c) \le B \}$. 
\begin{theorem}\label{thm_conv_randm_search}
\textit{For any $\delta>0$, define,
$$W(\delta):= \{\bPi \in \Gamma_c: |J_r(\bPi,\alpha_r)-J_r(\bPi^*,\alpha_r)| \le \delta\},$$
to be the set of $\delta$-optimal policies  for \eqref{eqn_original_problem}. Then $\{\bPi_k\}_{k\ge 1}$ generated from algorithm \ref{alg_random_search} visits $W(\delta)$ infinitely often with probability (w.p.) 1. Further $\hat{\bPi}^* \to W(\delta)$ a.s.} (see \eqref{eqn_pi_hat}).
\end{theorem}
\noindent {\bf Proof} is in \TR{\cite{TR}}{Appendix}. \eop

\section{Numerical Examples}
In this section, we present an example of inventory  control problem (see \cite{inv,puterman} for more details). There is an inventory with maximum possible size $M$. At the beginning of any day, the inventory owner needs to decide the additional quantity to be added to the  current inventory by ordering. If the owner orders, it pays a fixed ordering cost $O_f$ in addition to per unit cost $O_u$. We assume that any ordered quantity is received by the owner immediately.  The number of demands $D$ on any day is modelled by a geometric random variable with parameter $p$. Any unfulfilled demand is added towards the shortage cost, where per unit shortage cost is $C_s$. On the other hand, if there is any inventory remaining at the end of the day, the owner pays a holding cost $C_h$ per unit.  On any day $t$, the state of the system is given by current inventory level $x$, and  action $a$ represents the quantity to be ordered. The set of states is $\S_t = \S = \{0,\dots,M\}$, and set of  available actions $\A_{t,x}= \A_x$ (given $x$) equals  $\{0, \dots, M-x\}$.

In the first example, our aim is to minimize the risk-sensitive running objective-cost, that consists of ordering cost and holding cost, while keeping the risk-sensitive shortage cost below a bound $B$. The immediate running objective and constraint-cost are given by:
\begin{eqnarray*}
    r_t(x,a)  &=&  (O_f + a O_u)\indc{a>0} + C_h E[ (x +a - D)^+],\\
   c_t(x,a) &=&   C_s E[ (D- x - a )^+].
\end{eqnarray*}%
 Here,  we have replaced future state $x'$ dependent (through demand $D$) reward/constraint-cost $m(x,a,x')$ with expected reward/constraint-cost $m(x,a)$,   for simplicity by taking appropriate expectation. Nevertheless, this example provides good insights of the optimal policy for inventory control. We consider risk-sensitive framework as in \eqref{eqn_original_problem} to additionally minimize higher moments of the combined cost (as in \cite{uday})  with zero terminal objective/constraint-cost.
 
We study the variation in optimal objective-cost (value function) as the terminal time or risk-factor  varies.  We plot the normalized optimal objective-cost (value function) $v_r=\nicefrac{\log(J_r(\pi^*,\alpha_r))}{\gamma}$ and the constraint-cost $v_c=\nicefrac{\log(J_c(\pi^*,\alpha_c))}{\gamma_c}$; the other parameters are set to $M =5$, $O_f = 0.2$, $O_u = 0.4$, $C_h = 0.1$, $C_s =1$, $p=0.7$, $\beta=\beta_c = 0.8$, and initial distribution $\alpha_r(x) = \alpha_c(x) = \nicefrac{(M-x+1)}{\sum_{s}(M-s+1)}$. The bound on the risk-sensitive shortage cost is set to $B=e^{0.6 \gamma_c}$, with $\gamma_c = 0.1 \gamma$.
 In the left sub-figure of figure \ref{fig:performance}, the risk-factor for reward is set to $\gamma =0.5$.  One can see, when $T$ is small the problem is unconstrained. As $T$ increases, the problem becomes constrained, and then the   constraint-cost $v_c$ at optimality for any $T$ equals the (normalized) constraint $\nicefrac{\log(B)}{\gamma_c}$. Interestingly, the optimal objective initially increases with $T$, but then settles to a limit for higher values of $T$ as proved in \cite{uday}. 
 In the right sub-figure of figure \ref{fig:performance}, the terminal time is set to $T=5$.  We vary the risk-factor from 0.1 to 15, and observe that the reward decreases with $\gamma$.
  
 \TR{ \begin{figure}
    \centering
    \begin{minipage}{0.2 \textwidth}
       \includegraphics[trim={3cm 6cm 4cm 7cm},clip,scale=0.22]{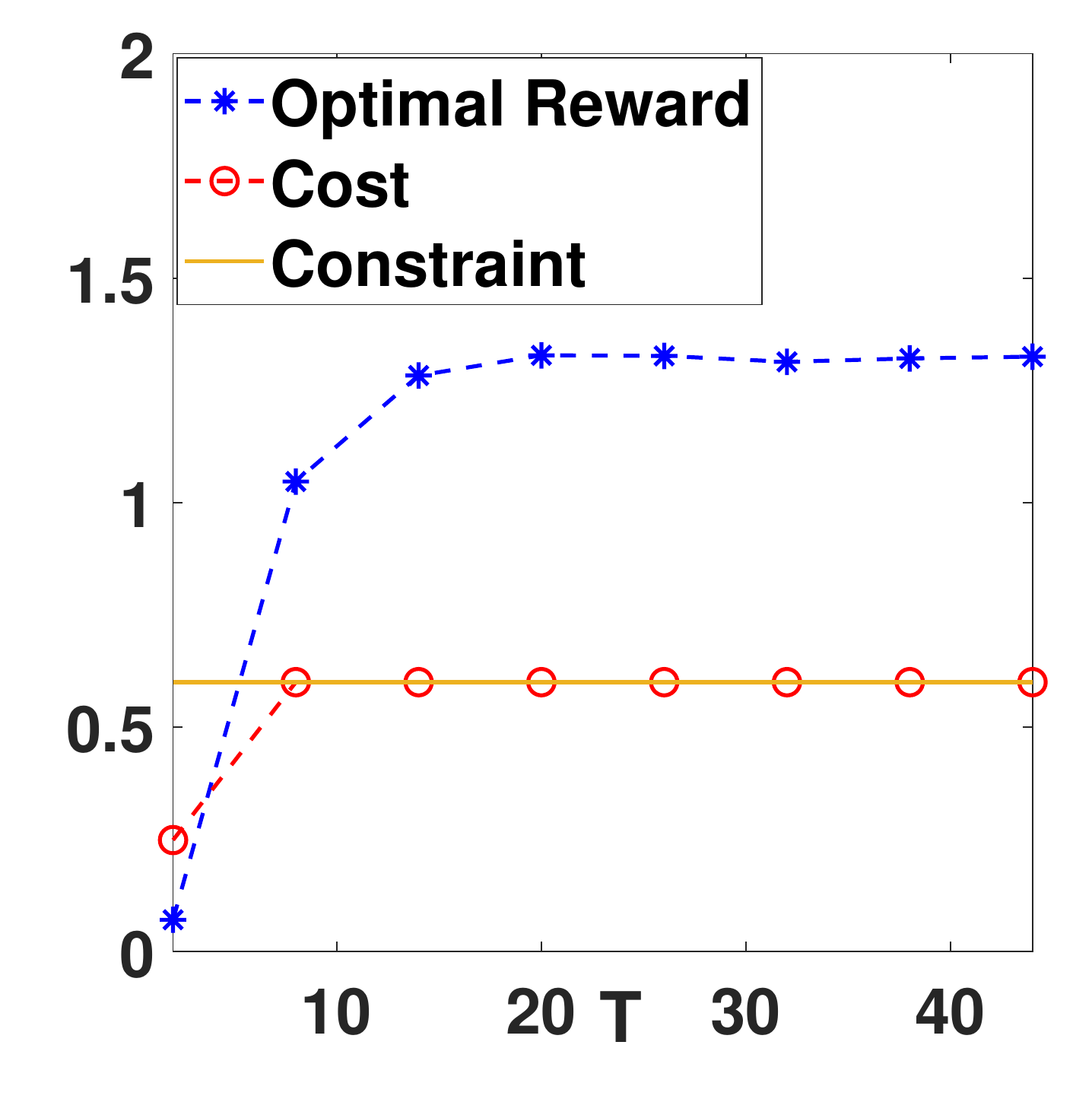}
    \end{minipage}%
    \hspace{2mm}
    \begin{minipage}{0.2 \textwidth}
    \includegraphics[trim={3cm 6cm 4cm 7cm},clip,scale=0.23]{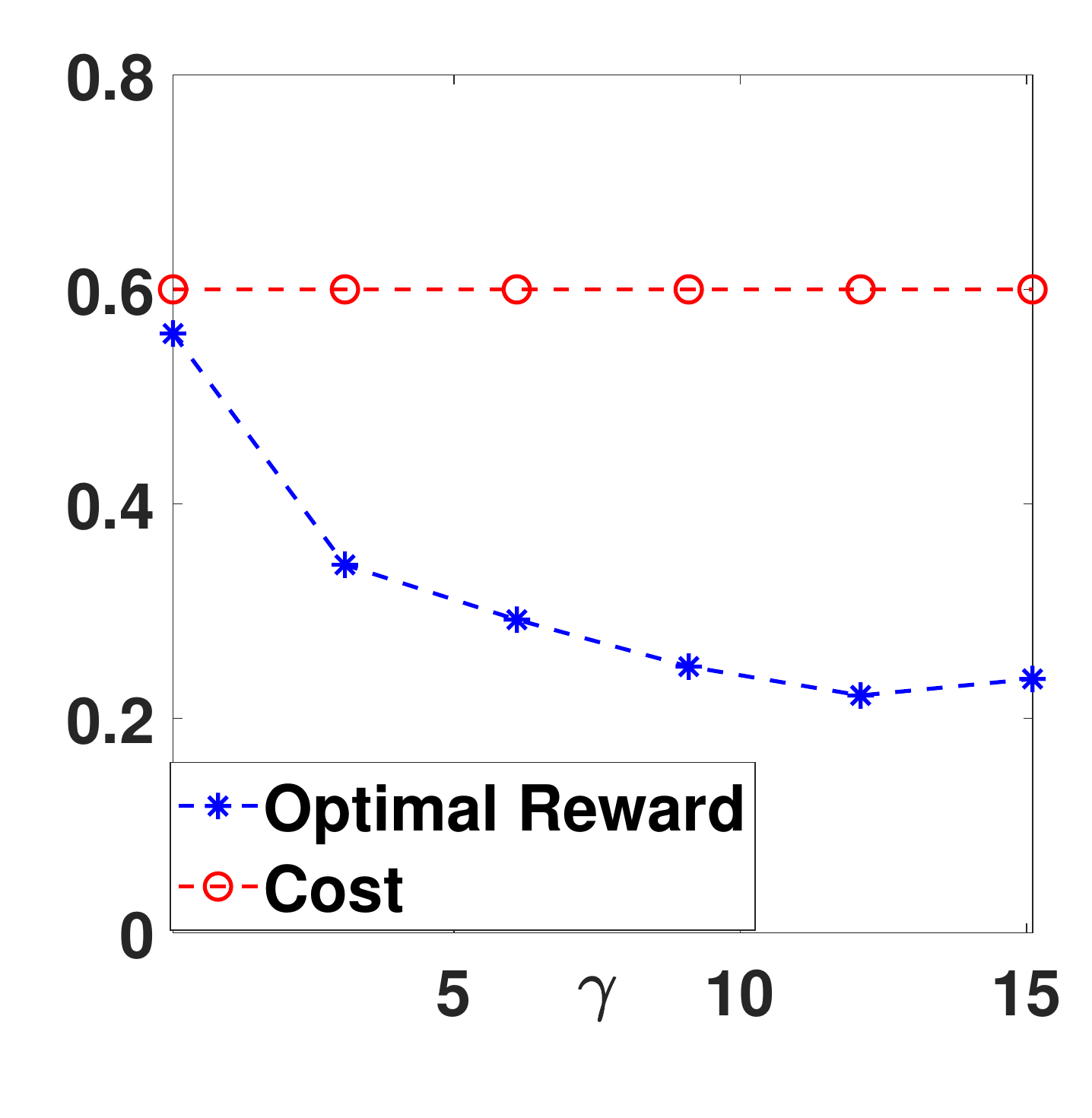}
        \end{minipage}
    \caption{Optimal running cost and shortage cost versus $T$ and risk-factor~$\gamma$}
    \label{fig:performance}
\end{figure}}%
{  \begin{figure}
    \centering
    \begin{minipage}{0.2 \textwidth}
       \includegraphics[trim={3cm 6cm 4cm 7cm},clip,scale=0.28]{performance_vs_T.pdf}
    \end{minipage}%
    \hspace{2mm}
    \begin{minipage}{0.2 \textwidth}
    \vspace{-1mm}
    \includegraphics[trim={3cm 6cm 4cm 7cm},clip,scale=0.31]{performance_vs_gam.pdf}
        \end{minipage}
    \caption{Optimal running cost and shortage cost versus $T$ and risk-factor~$\gamma$.}
    \label{fig:performance}
\end{figure}
}

We consider another example in inventory control with  running objective-cost consisting of ordering cost, holding cost and shortage cost, and a constraint on the number of orders.  The immediate running objective-cost and constraint-cost are now  given by:

\vspace{-3mm}
{\small
\begin{eqnarray*}
    r_t(x,a)  &=&  (O_f + a O_u)\indc{a>0}+ C_h E[ (x +a - D)^+] \\
    &&\hspace{1cm}  +C_s E[ (D- x - a )^+],\\
     c_t(x,a) &=& a.
\end{eqnarray*}}%
In figures \ref{fig:policy_with_gamma_pt5}-\ref{fig:policy_with_gamma_5}, we plot the optimal policy for different values of $\gamma$ and $B$. We   fix other parameters at $M =5$, $O_u = 0.2$, $C_s =6$, $p=0.6$, $\beta=\beta_c = 0.7$, $\gamma_c = \gamma$ with remaining parameters as in previous example.

In all the examples, the actions $0,1,2,3,4,5$ are represented respectively by blue circle, black star, red star, blue line, black line, red line.  For compact visualization of   the (non-stationary)  optimal policy we plot the  decision rule corresponding to state $x$  between values $[3(x+1), 3(x+1)+1]$; \textit{for example decision rule $d^*_t(x, a) = q$ is represented by a blue circle at $(t, 3(x+1)+q)$ when action $a = 0$ (ordering zero inventory).} 

\TR{\begin{figure*}
\begin{minipage}{0.31\textwidth}
    \begin{minipage}{0.4\textwidth}
      \includegraphics[trim={2cm 1cm 0cm 0.5cm},clip, scale=0.28]{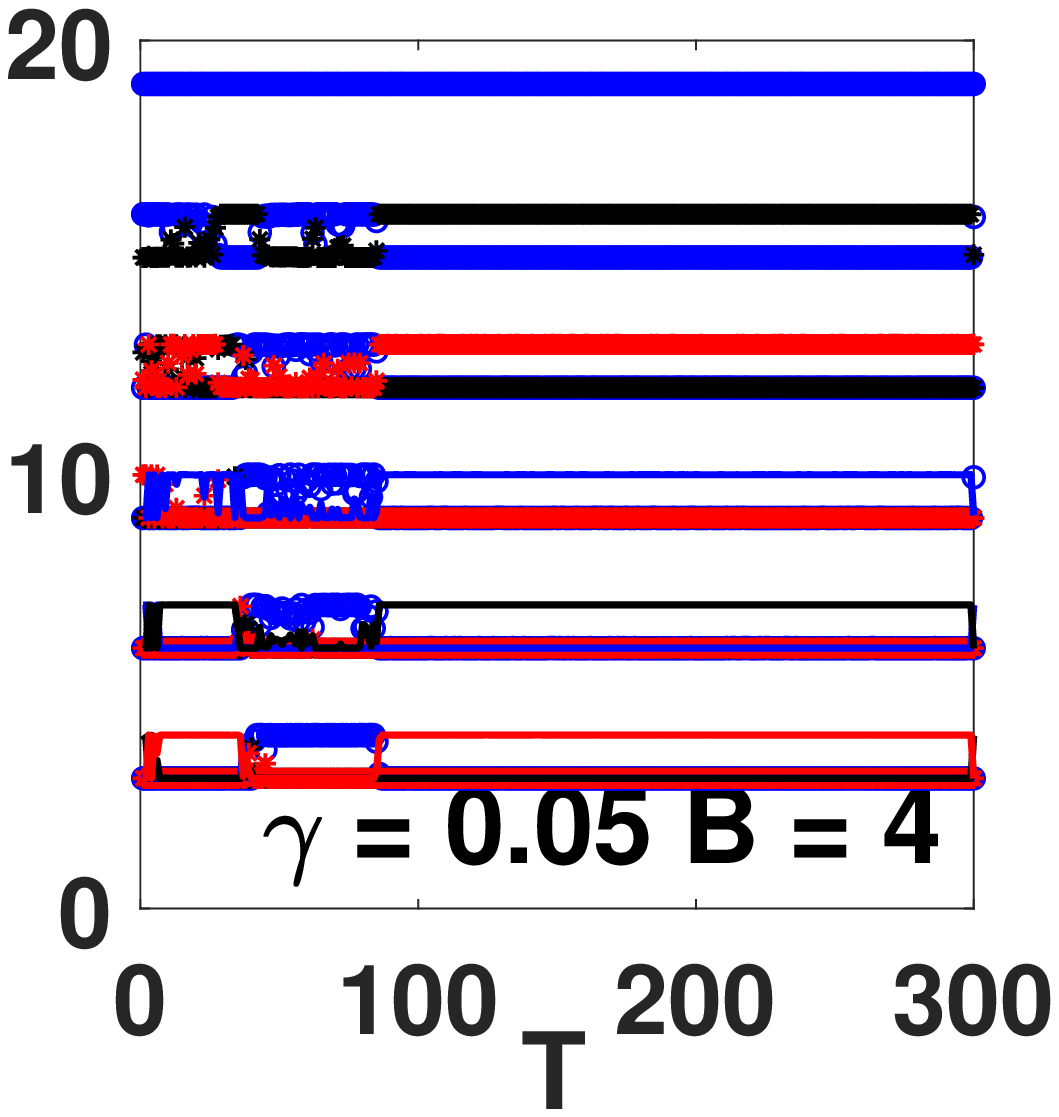}
    \end{minipage}\hspace{5mm}%
   \begin{minipage}{0.4 \textwidth}
   \vspace{-0.5mm}
     \includegraphics[trim={2.5cm 1cm 0cm 1cm},clip, scale=0.28]{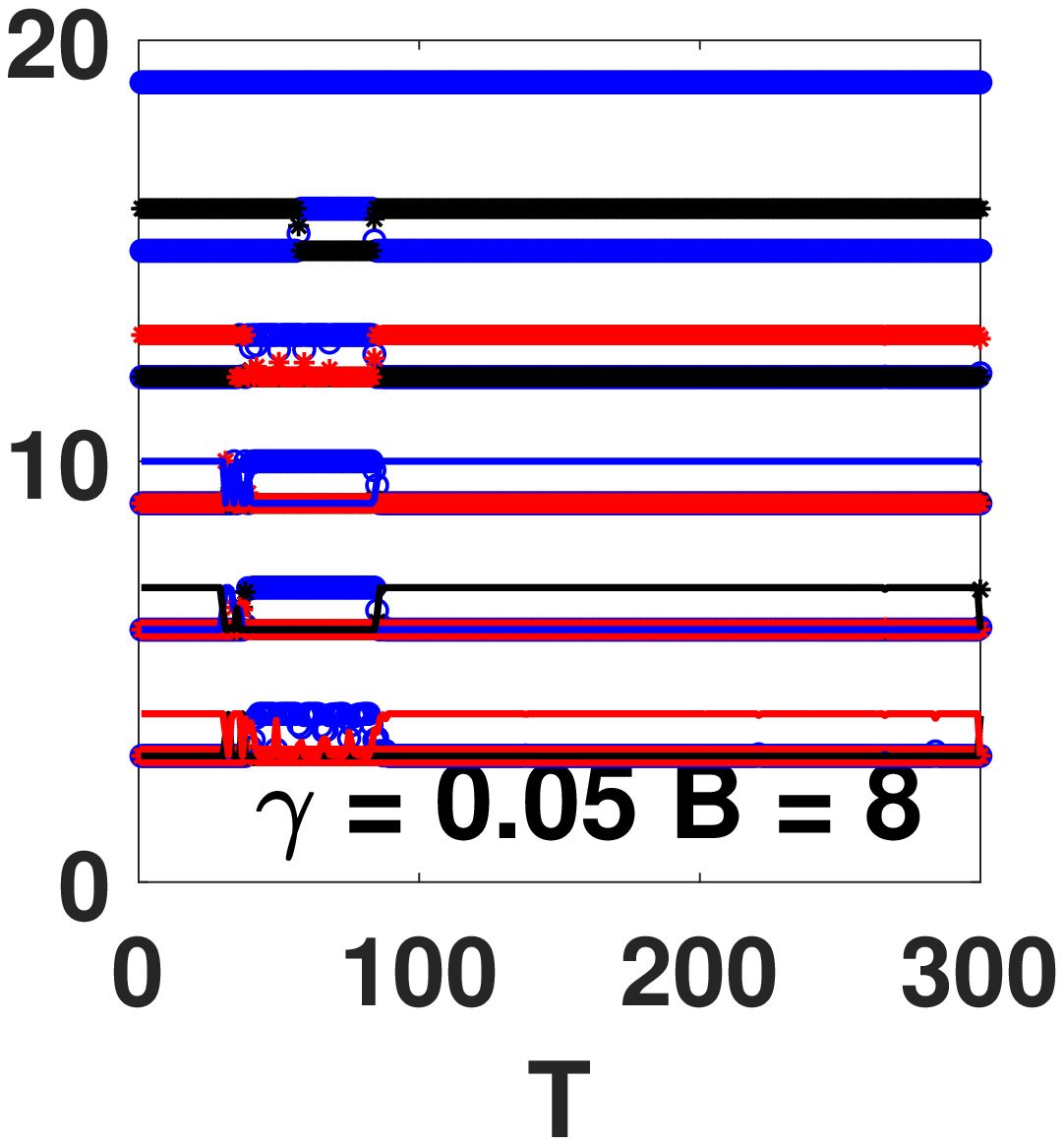}
    \end{minipage}
    \caption{Optimal policy with small risk-factor  \label{fig:policy_with_gamma_pt5}}%
\end{minipage}\hspace{2mm}%
\begin{minipage}{0.31\textwidth}
    \begin{minipage}{0.4 \textwidth}
\vspace{-0.5mm}
     \includegraphics[trim={2.8cm 1cm 0cm 0.2cm},clip, scale=0.28]{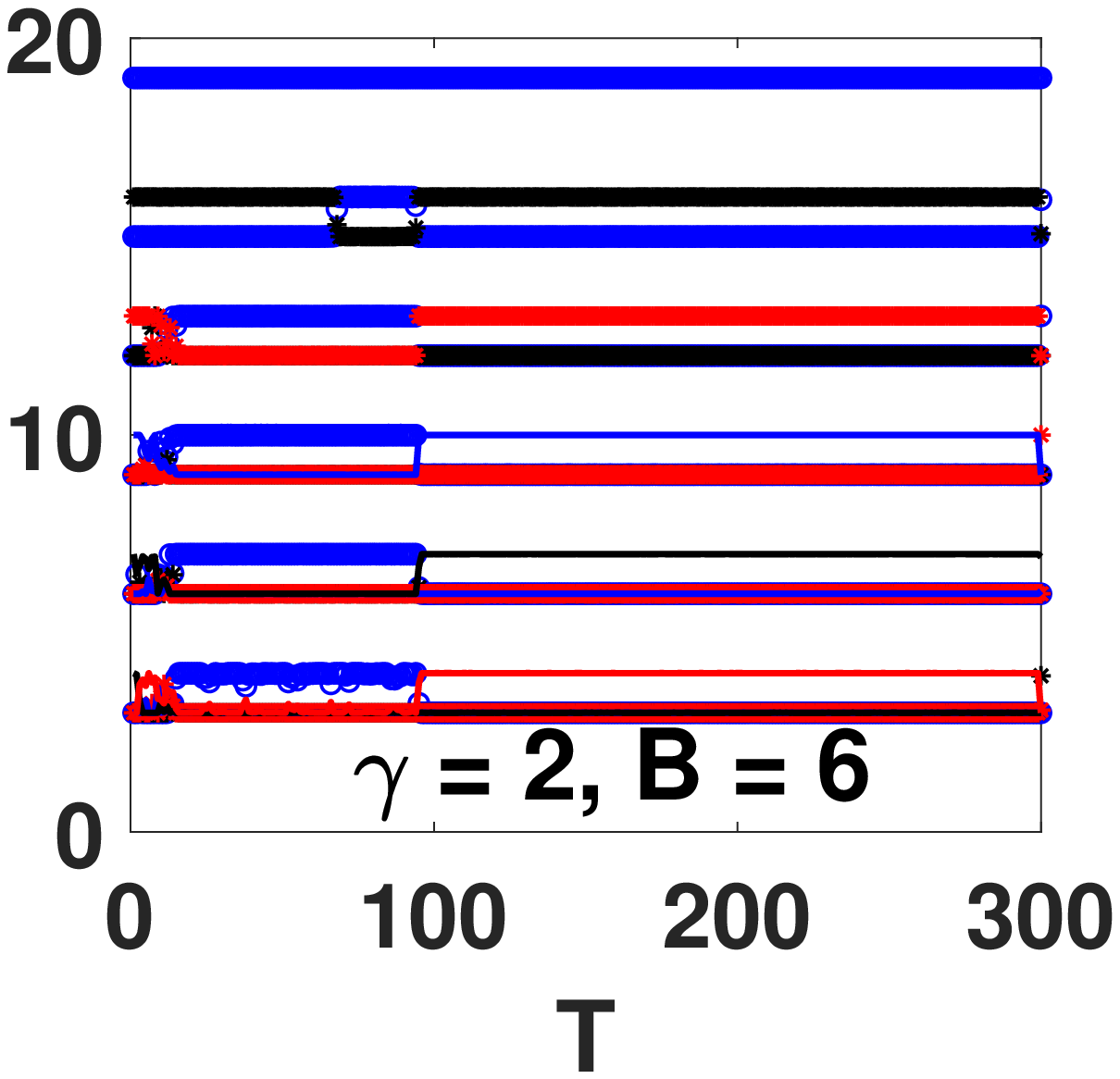}
    \end{minipage}\hspace{5mm}%
    \begin{minipage}{0.4\textwidth}
    \vspace{-2.3mm}
      \includegraphics[trim={3.5cm 1cm 0cm 0.6cm},clip, scale=0.28]{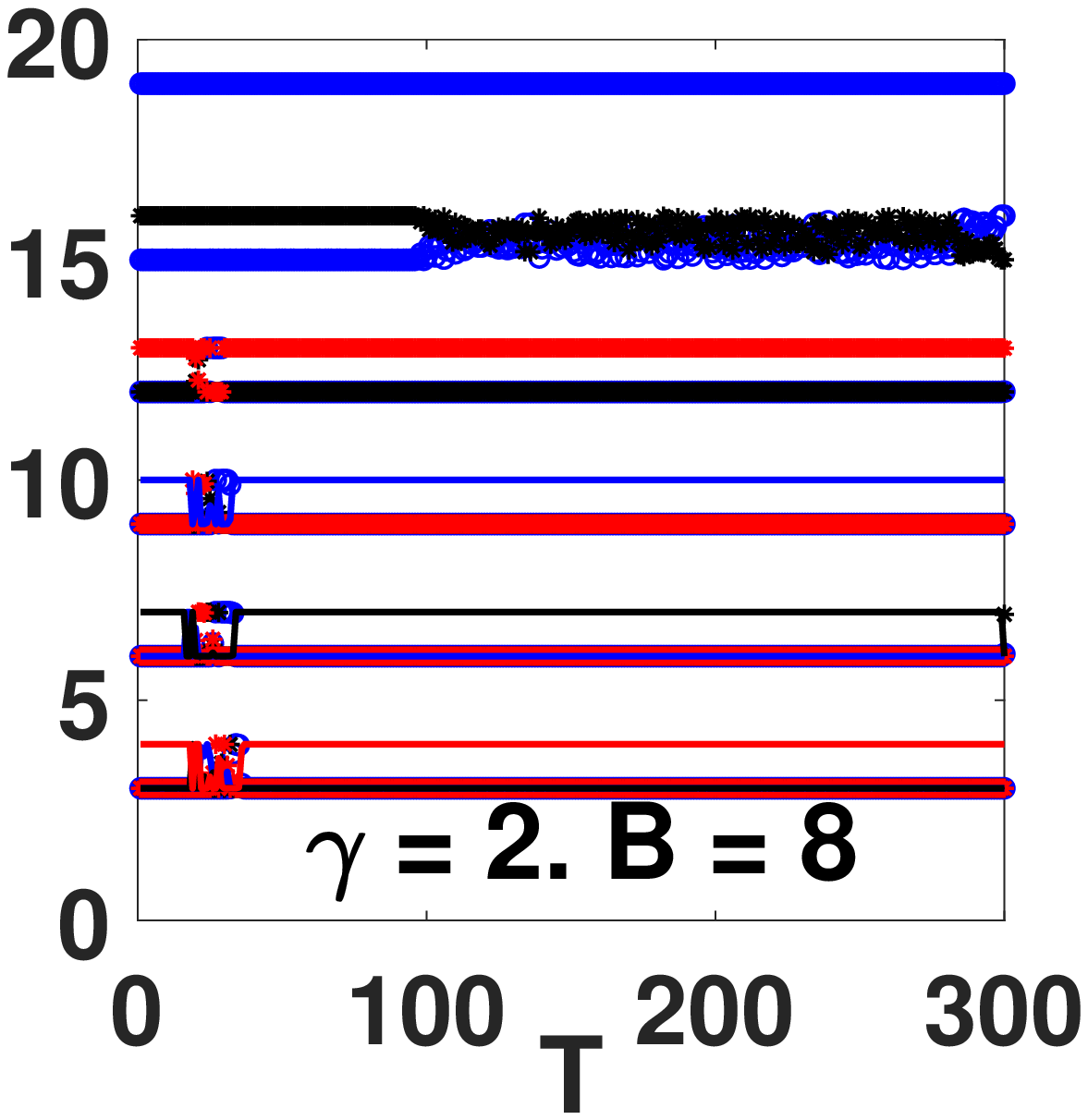}
    \end{minipage}
    \vspace{-0.8mm}
    \caption{With moderate risk-factor\label{fig:policy_with_gamma_2}}%
\end{minipage}\hspace{2mm}%
\begin{minipage}{0.31\textwidth}
\begin{minipage}{0.4\textwidth}
\vspace{-1mm}
      \includegraphics[trim={5.2cm 5.5cm 3cm 7.3cm},clip, scale=0.28]{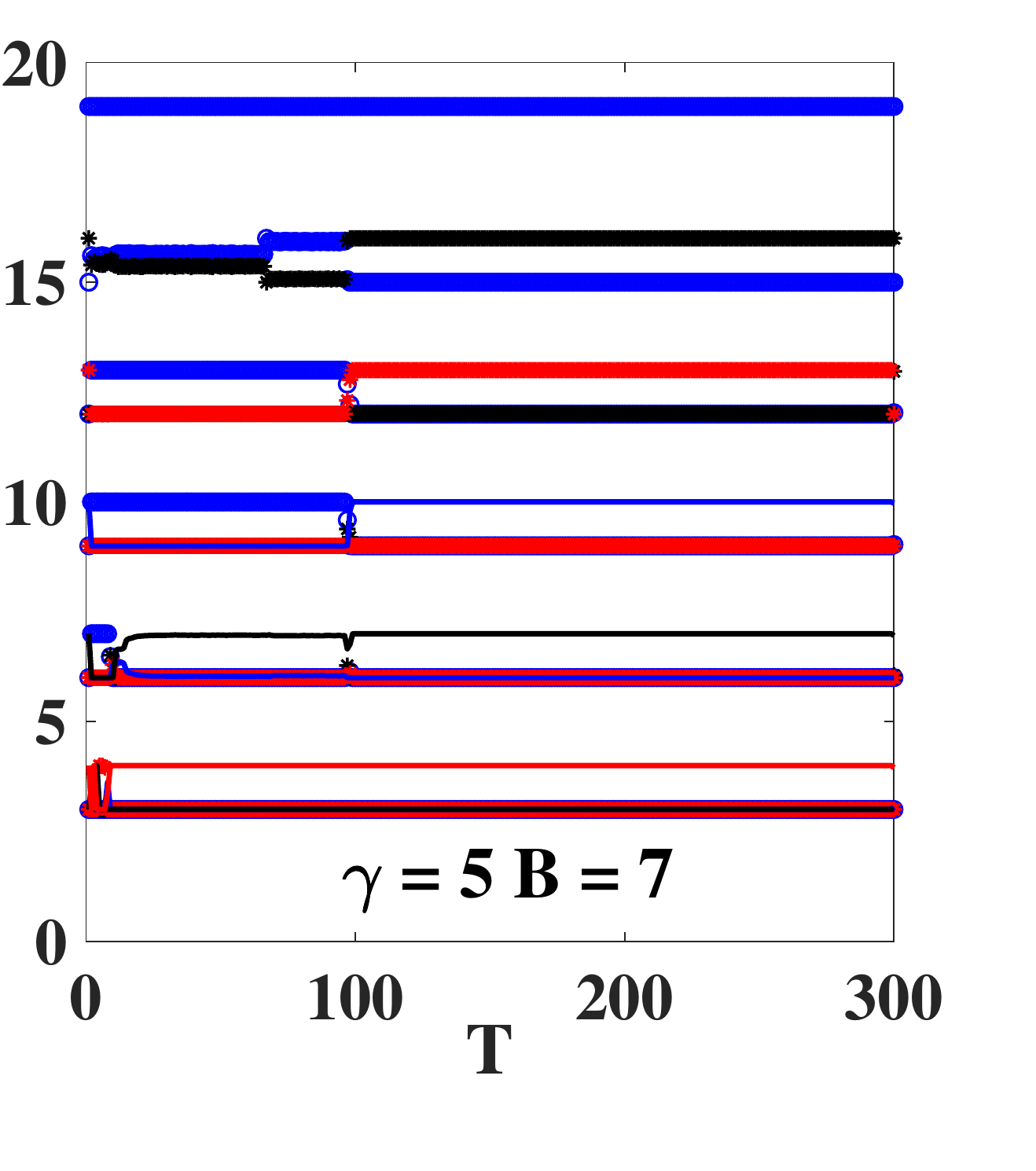}
    \end{minipage}\hspace{4mm}%
   \begin{minipage}{0.4 \textwidth}
   \vspace{0.8mm}
     \includegraphics[trim={3cm 5cm 4cm 6.5cm},clip, scale=0.28]{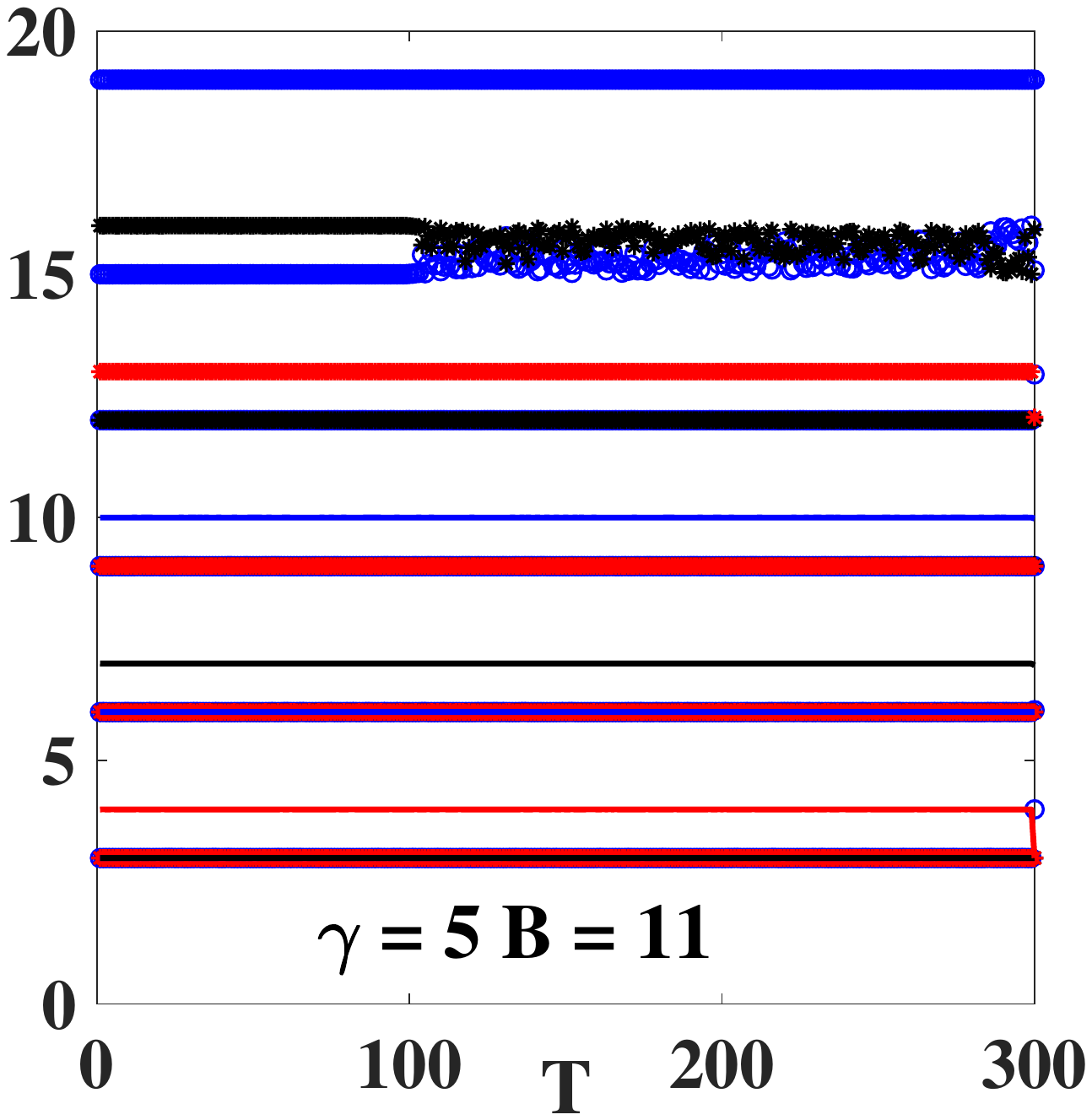}
    \end{minipage}
    \vspace{1mm}
    \caption{With high risk-factor \label{fig:policy_with_gamma_5}}
\end{minipage}
\end{figure*}
}%
{}

We see the optimal policy settles to a stationary rule  after initial non-stationary decision rules, i.e.,  \textit{we have an ultimately stationary policy in all the figures (as discussed in \cite{Jaquette76}). Interestingly this is true even for the problems with constraint and even the for risk factor as large as 5.} Further   more, we observe  ultimately stationary  deterministic policy (non-stationary and randomized only for initial time slots) for all the cases except for the right sub-figures of figures \ref{fig:policy_with_gamma_2},\ref{fig:policy_with_gamma_5}. 
These right sub-figures correspond to the case with larger risk-factor. 
It is well known that constrained problems are solved by randomized policies (e.g., \cite{eitan}); in all the left sub-figures (with active constraints) we see randomization in the initial time slots; further one orders less (more blue circles) with stricter constraints on the (risk-sensitive) number of orders. 
More interestingly with larger risk factor, the optimal policy orders less in higher states (compare right and left sub-figures for each case).  

We  plot this example with 300 time slots (in figures \ref{fig:policy_with_gamma_pt5}-\ref{fig:policy_with_gamma_5}), however we also obtained the policies with 1000 time-slots. The policy continues to  remain stationary after 300 time-slots. The pictures with 1000 time-slots lack clarity due to a large number of time slots and hence are avoided. 

These are   initial results and are provided just to illustrate the computational feasibility of the proposed algorithm. We obtain sufficiently good results within a few minutes for time slots even  as large as 100 time-slots. For 1000 time-slots the algorithm has to run for more than fifteen minutes to 
visualize a sufficiently accurate optimal policy. \textit{The last comment is that the  uniform distribution over corner points as ${\cal U}$ seems to perform better in many cases}.
\TR{}{\begin{figure}
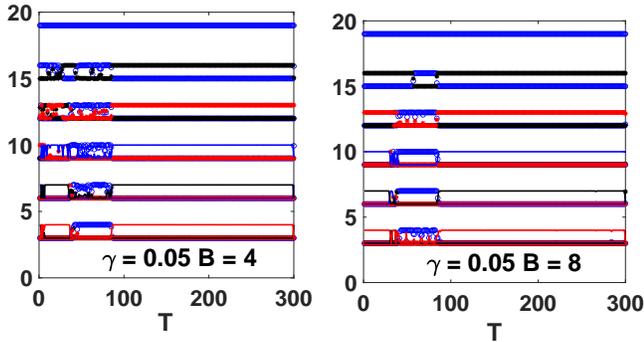
 
    \centering
    \begin{minipage}{0.23\textwidth}
      \includegraphics[trim={2cm 1cm 0cm 0.5cm},clip,scale=0.3]{policy_plots/T=300GamPt05_B_4_constrained.eps}
    \end{minipage}\hspace{2mm}%
   \begin{minipage}{0.23 \textwidth}
     \includegraphics[trim={2.5cm 1cm 0cm 1cm},clip,scale=0.3]{policy_plots/T=300GamPt05_B_8_unconstrained.eps}
    \end{minipage}
    \caption{Optimal policy with small risk-factors    \label{fig:policy_with_gamma_pt5}}
\end{figure}

\begin{figure}
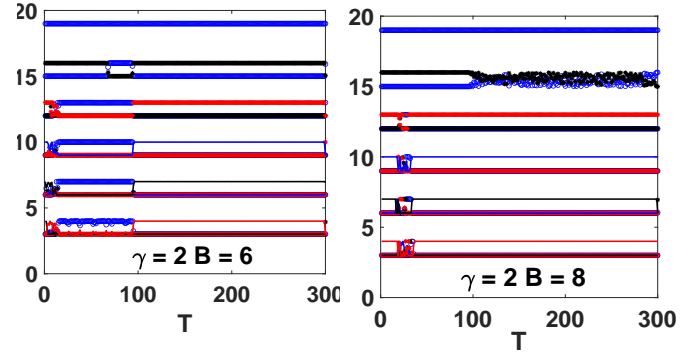
 
    \centering
    \begin{minipage}{0.23 \textwidth}
     \includegraphics[trim={3cm 0cm 0cm 0cm},clip,scale=0.3]{policy_plots/T=300Gam2_B_6.eps}
    \end{minipage}\hspace{2mm}%
    \begin{minipage}{0.23\textwidth}
      \includegraphics[trim={3.5cm 0cm 0cm 0cm},clip,scale=0.29]{policy_plots/T=300Gam2_B_8.eps}
    \end{minipage}
    \caption{Optimal policy with moderate risk-factor \label{fig:policy_with_gamma_2}}
\end{figure}

\begin{figure}
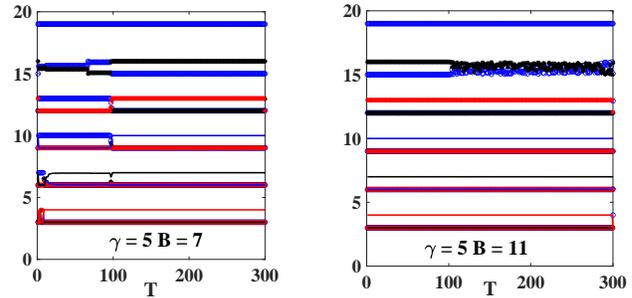
 
  \hspace{-5mm}  \begin{minipage}{0.23\textwidth}
      \includegraphics[trim={4cm 5.5cm 3cm 7.5cm},clip,scale=0.29]{policy_plots/T=300Gam5_B=7.pdf}
    \end{minipage}\hspace{2mm}%
   \begin{minipage}{0.23 \textwidth}
   \vspace{3mm}
     \includegraphics[trim={3cm 5cm 4cm 6.5cm},clip,scale=0.29]{policy_plots/T=300Gam5_B=11.pdf}
    \end{minipage}
    \caption{Optimal policy with high risk-factor \label{fig:policy_with_gamma_5}}
\end{figure}}

\section{Conclusions}
We consider a finite-horizon Risk-sensitive constrained MDP ($\mdp$),  and provide a fixed-point equation such that the optimal policy of $\mdp$ is also a solution. We further provide an iterative method to derive the solutions of the proposed fixed-point equation,  which involves solving a small Linear Programming (LP) problem at each iterate. We  propose a global optimization technique-based algorithm, namely GRC  algorithm, to derive the optimal policy for $\mdp$; this algorithm combines random restarts with fixed-point iterates of the proposed fixed-point equation using LP.  We show the convergence of the iterates in  the GRC algorithm to the optimal policy under certain conditions. The GRC algorithm is computationally feasible and the complexity grows only linearly with time-horizon of $\mdp$. We also present numerical results for risk-sensitive inventory control problems with a constraint.  We observe ultimately stationary policies are optimal (non-stationary decision rules    in initial time-slots and  a stationary   rule afterwards) for constrained problems in the presented examples.

This paper contains initial results to illustrate the idea to solve $\mdp$. One can strengthen the theoretical affirmation of the algorithm by proving few additional technical steps, which would be considered for the journal version of the paper.

We also attempted to understand the qualitative behaviour of the optimal policy for risk-sensitive, finite horizon and constrained inventory control problems. Some well known aspects like ultimately stationary policies, optimal $(S,s)$ policies, randomization with constraints etc., are observed. It would be   interesting to conduct a more elaborate study  of such properties, as now our algorithm made it possible to derive an optimal policy computationally.

\TR{}
{
\appendix
\section{Appendix}
\textbf{Proof of Lemma \ref{lem_exp_rew_pi_tilde}:} Let ${\mathcal{G}_t}:=({\bf X}_1^{t-1}, {\bf A}_1^{t-1})$, be the ordered tuple of random variables till time $(t-1)$. Further let $M_{\tau} :=m_\tau (X_\tau,A_\tau,X_{\tau+1})$ for $\tau\le T-1$, $M_T = m_T(X_T)$,  $\bPi'_t:=(\tilde{d}_t,
\d{t+1}{T-1})$ and $\bPi_t:=(\d{1}{t-1})$.
 From \eqref{eqn_f_lin},  for any $t\le T-1$, the LHS equals (using Markov property), 

\vspace{-3mm}
{\small
\begin{eqnarray*}
&=&  \sum_{x,a}  \theta^\bPi_{m,t}(x) \tilde{d}_t(x, a) Q^\bPi_{m,t}(x,a),\\
 &=& \sum_{x,a}  \theta^\bPi_{m,t}(x) \tilde{d}_t(x, a)  E_{\alpha_m}^\bPi\left[e^{\sum_{\tau=t}^{T}M_\tau}| X_t=x,A_t=a \right],\\
 &=& \sum_{x}  \theta^\bPi_{m,t}(x)  E^{\bPi'_t}_{\alpha_m}\left[e^{\sum_{\tau=t}^{T}M_\tau}| X_t=x \right],\\
&=& \sum_{x}  E_{\alpha_m}^{\bPi_t}\left[e^{\sum_{\tau=1}^{t-1}M_\tau} \indc{X_t =x} E_{\alpha_m}^{\bPi'_t}\left[e^{\sum_{\tau=t}^{T}M_\tau}|X_t=x \right]\right], \\
&=&   E_{\alpha_m}^{\bPi_t}\left[\sum_{x} \indc{X_t =x} E_{\alpha_m}^{\bPi'_t}\left[e^{\sum_{\tau=1}^{T}M_\tau}| \mathcal{G}_t, X_t=x \right]\right], \\
&=&   E_{\alpha_m}^{\bPi_t}\left[ E_{\alpha_m}^{\bPi'_t}\left[e^{\sum_{\tau=1}^{T}M_\tau}| \mathcal{G}_t , X_t\right]\right],\\
&=&   E_{\alpha_m}^{\etapi}\left[ E_{\alpha_m}^{\etapi}\left[e^{\sum_{\tau=1}^{T}M_\tau}| \mathcal{G}_t,X_t \right]\right]=  E_{\alpha_m}^{\etapi}\left[e^{\sum_{\tau=1}^{T}M_\tau}\right].
\end{eqnarray*}}Hence proved.  \eop


\textbf{Proof of Theorem \ref{thm_nec_cndi}:}
Let $\bPi^*$ be the optimal policy for $\mdp$ \eqref{eqn_original_problem} and $\bbtheta^*_m, \bbQ^*_m$ be the corresponding forward and backward factors. Consider $\lp{\bPi^*}$, the feasible region is non-empty because $\bPi^*$ is feasible  from Lemma \ref{lem_exp_rew_pi_tilde} and \eqref{eqn_original_problem}. Consider any $\bPi \in \M(\bPi^*)$, then for any $t\le T-1$,
\begin{eqnarray*}
f_t(\bPi, \bbtheta_r^{*}, \bbQ^*_r) &=& J_r(\eta^t_{\bPi^*,\bPi}, \alpha_r),\\
&\le& J_r(\bPi^*, \alpha_r)= f_t(\bPi^*, \bbtheta_r^{*}, \bbQ^*_r).
\end{eqnarray*}%
Hence
$\sum_{t=1}^{T-1} f_t(\bPi, \bbtheta_r^*, \bbQ^*_r)  \le\sum_{t=1}^{T-1} f_t(\bPi^*, \bbtheta_r^*, \bbQ^*_r)$, implying $\bPi^* \in \M(\bPi^*)$\eop

\textbf{Proof of Theorem \ref{thm_eqvlnt}:}
The proof is in two steps. First we  show the equivalency of $\mdp$ and \textbf{GF} problem. 

Let $\bPi^*$ be the optimal policy of $\mdp$ \eqref{eqn_original_problem}. From Theorem \ref{thm_nec_cndi}, $\bPi^*\in \Mset$. Thus,
$$\sup_{\Pi \in \Mset}\sum_{t=1}^{T-1} f_t({\bPi}, \bbtheta_r^{\bPi}, \bbQ_r^{\bPi}) \ge \sum_{t=1}^{T-1} f_t({\bPi^*}, \bbtheta_r^{*}, \bbQ_r^{*}).$$
Any policy $\bPi \in \Mset$ satisfies $\bPi\in\M(\bPi)$ by definition of $\Mset$. By lemma \ref{lem_exp_rew_pi_tilde}, $J_c(\bPi,\alpha_c) = f_t(\bPi, \bbtheta_c^{\bPi}, \bbQ_c^{\bPi})\le B$ from the constraint in $\lp{\bPi}$. So, $\bPi$ is feasible for $\mdp$ \eqref{eqn_original_problem}. Thus, from optimality of $\Pi^*$ and
Lemma \ref{lem_exp_rew_pi_tilde},
$$\sup_{\Pi \in \Mset}\sum_{t=1}^{T-1} f_t({\bPi}, \bbtheta_r^{\bPi}, \bbQ_r^{\bPi}) \le \sum_{t=1}^{T-1} f_t({\bPi^*}, \bbtheta_r^{*}, \bbQ_r^{*}).$$
Thus, the equality holds and we have the result.

Towards the equivalency of $\mdp$ \eqref{eqn_original_problem} and \textbf{GO} problem, observe that any feasible policy $\bPi$ for the \textbf{GO}  problem is also a feasible policy for \eqref{eqn_original_problem} and vice-versa from Lemma \ref{lem_exp_rew_pi_tilde}. Hence, we have the results. \eop

\textbf{Proof of Theorem \ref{thm_cnvg_ode}:}
To prove this theorem, we will first show that the 
following sequence  of piece-wise constant functions that  start with $\bPi_k$  are equicontinuous in an extended sense (referred as \textit{equicontinuous} for brevity):  

\vspace{-3mm}
 {\small \begin{eqnarray*}
         \bPi^k(t) :=  \bPi_k +\hspace{-2mm} \sum_{j=k}^{\nu(t_k+t)-1}\hspace{-5mm} \epsilon_j (\rpi_j - \bPi_j),  
\end{eqnarray*}}with $\nu(t):=\max\{k:t_k\le t\}$ and $\rpi_j$ is the policy randomly chosen from $\M(\bPi_{j})$ (see \eqref{eqn_local_update}). Then the result follows from  \cite[Chapter 5, Theorem 2.2]{kushner2003stochastic}.  
Recall $\bPi^k(t)$ is a vector (one component for each state $x$, action $a$ and decision epoch $\tau$), and let $i$ be one such (arbitrary) component  represented briefly by
 $  \cPii^k(t)$. The function  $  \cPii^k(t)$, then equals (see \eqref{eqn_local_update}):

   \vspace{-3mm}
 {\small \begin{eqnarray*}
          \cPii^k(t) :=  \dPii_k +\hspace{-2mm} \sum_{j=k}^{\nu(t_k+t)-1}\hspace{-5mm} \epsilon_j (\Psii_j - \dPii_j),  
\end{eqnarray*}}with $\Psii_j$ representing the $i$-th component of $\rpi_j$. 
We will prove that every component of function $\bPi^k(\cdot)$ is equicontinuous. Basically, we need to show that, for each $T$ and $\varepsilon>0$, there is a $\delta>0$ such that,
\begin{eqnarray*}
             \lim\sup_k\sup_{0\le t-s\le \delta,|t|<T }|\cPii^k(t)-\cPii^k(s)|<\varepsilon \mbox{ for every }i.
        \end{eqnarray*}
   The proof of equicontinuity is exactly similar to that provided in the proof of \cite[Chapter 5, Theorem 2.1]{kushner2003stochastic} for the case with continuous $g$, except for the fact that $g (\cdot)$ in our case is not continuous (see \eqref{eqn_ode}). We will only provide differences in the proof steps towards 
       $\{\cPii^k(t)\}_k$ sequence, and it can be proved analogously for others.
       
Define $\delta M_k := \Psii_k - \dPii_k - \gi(\bPi_k) $, where $\gi(\cdot)$ is the $i$-th component of function $g(\cdot)$ of \eqref{eqn_ode}. Then function $\cPii^k(t)$ can be re-written as,
     \begin{eqnarray*}
          \cPii^k(t)=\dPii_k + \sum_{j=k}^{\nu(t_k+t)-1} \epsilon_j \delta M_j+ \sum_{j=k}^{\nu(t_k+t)-1}\epsilon_j \gi (\bPi_j),
     \end{eqnarray*}     Now define $M_k=\sum_{\tau=0}^{k-1}\epsilon_\tau\delta M_\tau $, it is easy to prove $\{M_k,\mathcal{F}_k\}$ is martingale, where $\mathcal{F}_k$ is natural filtration. Thus, using Martingale inequality (see \cite[Chapter 4, equation (1.4)]{kushner2003stochastic} for $q(M)=M^2$), we get for each $\mu>0$,
       \begin{equation*}
           P_{\mathcal{F}_\tau}\left\{\sup_{\tau\le j \le k}|M_j-M_\tau|\ge \mu\right\}\le \frac{E_{\mathcal{F}_\tau}|\sum_{m=\tau}^{k-1}\epsilon_m\delta M_m|^2}{\mu^2}.
       \end{equation*} Using the fact $E[\delta M_\tau \delta M_j
       ]=0$ for $\tau<j$,  and that     $\sup_k E[\delta M_k^2] \le c$ for some constant $c$, since all the involved quantities (policies) are upper-bounded by 1, we have,
           \begin{eqnarray} \label{eqn_mtg_ineq}
  \lim_{\tau\rightarrow\infty}  P\left\{\sup_{\tau\le j }|M_j-M_\tau| \geq \mu \right\}=0 \mbox{ for each } \mu>0.
               \end{eqnarray}   
        Further, we can re-write $\cPii^k(t)$ as,
        \begin{equation*}
            \cPii^k(t)=\dPii_k+\int_0^t \gi(\bPi^k(z)) dz + M^{k} (t) + \rho^{k} (t) 
        \end{equation*}
        where we denote $ M^{k} (t)= \sum_{\tau=k}^{\nu(t_k+t)-1} \epsilon_\tau\delta M_\tau $, and $ \rho^{k} (t) =\sum_{\tau=k}^{\nu(t_k+t)-1} \epsilon_\tau \gi(\bPi_\tau)-\int_0^t \gi(\bPi^k(z)) dz$. 
        
        Clearly $\cPii^k(0)=\dPii_k \le 1$ and thus to claim equi-continuity, we have,
        
        \vspace{-3mm}
  {\small \begin{eqnarray} \label{eqn_terms}
     \sup  |\cPii^k(t)-\cPii^k(s)|&\le&\sup\left|\int_s^t \gi(\bPi^k(z)) dz\right| \nonumber\\   
   && +\sup\left| M^{k} (t) - M^{k} (s) \right|\nonumber\\   
   && + \sup\left| \rho^{k} (t)- \rho^{k} (s)\right|,
   \end{eqnarray}}%
where supremum is taken over $S_T := \{(s,t): 0\le t-s\le \delta,|s|<T,|t|<T\}$. Let us consider the first term from above, since $\gi(\bPi^k(z))\le 1$,
\begin{eqnarray*}
     \left|\int_s^t \gi(\bPi^k(z)) dz \right|  
       \le 2(t-s)\le \delta. 
\end{eqnarray*}
Towards the second term, first by \eqref{eqn_mtg_ineq} and continuity of probability,
\begin{eqnarray*}
        P\left\{\lim_{\tau\rightarrow\infty}\sup_{\tau\le j }|M_j-M_\tau|\ge\mu\right\}&=&0 \mbox{ for all } \mu>0.
\end{eqnarray*}
Let $A_n:=\left\{\omega: \lim_{\tau\rightarrow\infty}\sup_{\tau\le j }|M_j-M_\tau|<\frac{1}{n}\right\}$, then $P(A_n)=1$ for each $n>0$. Now we claim that for any 
$$\omega\in \cap_{n>0}A_n \mbox{ and any } T,\ \ \  \sup_{t \le T} M^{k}(t)\rightarrow 0,$$
as then the second term in \eqref{eqn_terms} converges to zero because  $ |M^{k} (t) - M^{k} (s)| \le  |M^{k} (t) | + |M^{k} (s)|$.
To this end, for every $\omega\in \cap_{n>0}A_n$,
\begin{eqnarray*}
     \sup_{t  \le T}|M^{k}(t)|&=& \sup_{t \le T}|M_{\nu(t_k+t)}-M_k| \le \sup_{j \ge k} |M_{j}-M_k|.
\end{eqnarray*}as $\nu(t_k+t) \ge k$ in the above. By taking $k \to \infty$ first, the above is upper bounded by $1/n$  (see definition of $A_n$) for each $n$. Then letting $n\to\infty$, we get our claim. 

 For the last term, it can be proved by induction that when $t$ exactly corresponds to the end of epochs, i.e., when $t = t_n-t_k, \ (n>k)$ that $\rho^{k}(t)=0$.  Now, we are only left to prove that $\rho^{k}(t)\rightarrow 0$ uniformly in $t$ (for general $t$) as $k \rightarrow \infty$. We will prove this claim for each $T>0$, such that $|t|<T$:

\vspace{-2mm}
{\small
\begin{eqnarray*}
 |\rho^{k}(t)| &=&\left| \sum_{\tau=k}^{\nu(t_k+t)-1} \epsilon_\tau \gi(\bPi_\tau)-\int_0^t \gi(\bPi^k(z)) dz\right|,\\
 &=&\left|\int_{\nu(t_k+t)}^t   \gi(\bPi^k(z)) dz\right|, \\
 &\le& (t - n(t_k + t))\ \leq \epsilon_k,
\end{eqnarray*}}where $\epsilon_k \to 0$ and first inequality follows as in proof of first term. This proves the equicontinuity in extended sense for $\{\cPii^k(t)\}_k$. Proof follows in exact similar lines for all components. This proves $\{\bPi^k(t)\}_k$ is equicontinuous in extended sense. Hence, from \cite[Chapter 5, Theorem 2.2]{kushner2003stochastic} $\bPi_k\rightarrow \mathbb{A}$.
\eop

\textbf{Proof of Theorem \ref{thm_conv_randm_search}:} It is sufficient to prove that the sequence of policies generated in Algorithm \ref{alg_random_search} visits the set $W(\delta)$ for any $\delta>0$. Towards this, we apply the results in \cite[Theorem 1]{global}, which says the sequence $\{\bPi_j\}$ visits the set $W(\delta)$ infinitely often with probability 1 if the following conditions are satisfied: i) $\Gamma_c$ (feasible region) is a compact set, ii) $J_r(\cdot,\alpha_r)$ is Lipschitz continuous, iii) and $\sum_{j=1}^{\infty} q_j(\varepsilon) =\infty$ for any $\varepsilon>0$, where $q_j(\varepsilon) := \inf_\bPi P_j(B(\bPi,\varepsilon))$ with  $P_j(B(\bPi,\varepsilon))$ as the probability that $\Pi_j$ (in algorithm \ref{alg_random_search}) is chosen from  $\varepsilon$-ball centered at $\bPi$, $B(\bPi,\varepsilon)=\{\bPi'\in \Gamma_c : \parallel\bPi-\bPi' \parallel\le \varepsilon\}$. The first two conditions are trivially true because of finite action and state space. Further, since we generate new policy uniformly from the whole space with probability $p_k:= \frac{w}{k}$ for a constant $w$, the third condition is also true. Hence we have the proof. \eop }
\end{document}